
\magnification\magstep1
\hoffset=0.5truecm
\voffset=0.5truecm
\hsize=15.8truecm
\vsize=23.truecm
\baselineskip=14pt plus0.1pt minus0.1pt \parindent=19pt
\lineskip=4pt\lineskiplimit=0.1pt      \parskip=0.1pt plus1pt

\font\ftitle=cmbx10 scaled\magstep1
\font\fabs=cmbx10 scaled800

\hfill December 1992
\vskip2truecm
\centerline{\ftitle MODULAR STRUCTURE AND}\smallskip
\centerline{\ftitle DUALITY IN CONFORMAL}\smallskip
\centerline{\ftitle QUANTUM FIELD THEORY}\bigskip\bigskip
\centerline{R. Brunetti$^1$\footnote{$^*$}{
Supported in part by Ministero della Ricerca Scientifica and CNR-GNAFA. }
\footnote{$^\bullet$}{Supported in part by INFN, sez. Napoli.},
D. Guido$^{2*}$ and R. Longo$^{2*}$}
\vskip1.truecm
\item{$(^1)$}  Dipartimento di Fisica, Universit\`a di Napoli ``Federico II''
\par Mostra d'Oltremare, Pad. 19\qquad I--80125, Napoli, Italy
\par E-mail BRUNETTI@NAPOLI.INFN.IT
\item{$(^2)$} Dipartimento di Matematica, Universit\`a di Roma ``Tor Vergata''
\par Via della Ricerca Scientifica\qquad I--00133, Roma, Italy.
\par E-mail GUIDO@MAT.UTOVRM.IT,\quad LONGO@MAT.UTOVRM.IT
\vskip1.truecm

\noindent {\fabs ABSTRACT. Making use of a recent result of Borchers, an
algebraic version of the Bisognano-Wichmann  theorem is given for  conformal
quantum field theories, i.e. the Tomita-Takesaki modular group
associated with the von Neumann algebra of a wedge region and the
vacuum vector concides with the evolution given by the rescaled pure
Lorentz transformations
preserving the wedge. A similar geometric description is valid
for the algebras associated with double cones. Moreover
essential duality holds  on the Minkowski space $M$, and Haag
duality for double cones holds provided  the net of local algebras
is extended to a pre-cosheaf on the superworld $\tilde M$, i.e. the
universal covering of the Dirac-Weyl compactification of $M$. As a
consequence a PCT symmetry exists for any algebraic conformal field
theory in even space-time dimension. Analogous results hold for a
Poincar\'e covariant theory provided the modular groups
corresponding to wedge algebras have the expected geometrical
meaning and the split
property is satisfied. In particular the Poincar\'e representation
is unique in this case.}

\vfill\eject

\def\A{{\cal A}}
\def\tildA{\tilde{\cal A}}
\def\tildAnul{\tilde{\cal A}_0}
\def\a{\alpha}
\def\B{{\cal B}}
\def\b{\beta}
\def\C{{\cal C}}
\def\Co{{\bf C}}
\def\D{\Delta}
\def\E{{\cal E}}
\def\f{\varphi}
\def\g{\gamma}
\def\G{{\cal G}}
\def\H{{\cal H}}
\def\K{{\cal K}}
\def\L{\Lambda}
\def\m{\mu}
\def\Na{{\bf N}}

\def\O{{\cal O}}
\def\p{\pi}
\def\Q{\Omega}
\def\r{\rho}
\def\R{{\cal R}}
\def\Re{{\bf R}}

\def\t{\tau}
\def\T{\Theta}

\def\W{{\cal W}}
\def\x{\xi}

\def\Ze{{\bf Z}}
\def\ad{{\rm ad}}
\def\Lor{{\cal L}_+^\uparrow}
\def\Poi{{\cal P}_+^\uparrow}

\def\imply{\Rightarrow}
\def\np{\par\noindent}
\catcode`@=11
\def\quad@rato#1#2{{\vcenter{\vbox{
        \hrule height#2pt
        \hbox{\vrule width#2pt height#1pt \kern#1pt \vrule width#2pt}
        \hrule height#2pt} }}}
\def\quadratello{\mathchoice
\quad@rato5{.5}\quad@rato5{.5}\quad@rato{3.5}{.35}\quad@rato{2.5}{.25} }
\catcode`@=12
\def\proof#1\par{\medskip\noindent{\bf Proof}\quad#1\hfill$\quadratello$
    \bigskip}
\def\proofof#1#2\par{\medskip\noindent
   {\bf Proof of #1}\quad#2\hfill$\quadratello$\bigskip}
\def\section #1\par{\vskip0pt plus.3\vsize\penalty-75
    \vskip0pt plus -.3\vsize\bigskip\bigskip
    \noindent{\sectionfont #1}\nobreak\smallskip\noindent}
\def\subs#1{\medskip\noindent{\it#1}\qquad}
\def\claim#1#2\par{\vskip.1in\medbreak\noindent{\bf #1.} {\sl #2}\par
    \ifdim\lastskip<\medskipamount\removelastskip\penalty55\medskip\fi}
\def\rmclaim#1#2\par{\vskip.1in\medbreak\noindent{\bf #1.} { #2}\par
    \ifdim\lastskip<\medskipamount\removelastskip\penalty55\medskip\fi}
\font\sectionfont=cmbx10 scaled\magstep1
\newcount\REFcount \REFcount=1
\def\numref{\number\REFcount}
\def\addref{\global\advance\REFcount by 1}
\def\wdef#1#2{\expandafter\xdef\csname#1\endcsname{#2}}
\def\wdch#1#2#3{\ifundef{#1#2}\wdef{#1#2}{#3}
    \else\write16{!!doubly defined#1,#2}\fi}
\def\wval#1{\csname#1\endcsname}
\def\ifundef#1{\expandafter\ifx\csname#1\endcsname\relax}
\def\ref(#1){\wdef{q#1}{yes}\ifundef{r#1}$\diamondsuit$#1
  \write16{!!ref #1 was never defined!!}\else\wval{r#1}\fi}
\def\inputreferences{
    \def\REF(##1)##2\endREF{\wdch{r}{##1}{\numref}\addref}
    \REFERENCES}
\def\references{
    \def\REF(##1)##2\endREF{\ifundef{q##1}
        \write16{!!ref. [##1] was never quoted!!}\fi
        \item{[\ref(##1)]}##2}
    \section{References}\par\REFERENCES}
\def\REFERENCES{
\REF(Arak1)Araki H. , ``{\it A lattice of von Neumann
algebras associated with the quantum field theory of a
free Bose field}", J. Math. Phys. 4 (1963), 1343-1362.
\endREF
\REF(BiWi1)Bisognano J., Wichmann E. , ``{\it On the duality condition
for a Hermitian scalar field}", J. Math. Phys. 16 (1975), 985-1007
\endREF
\REF(Borc1)Borchers H.J. , ``{\it The CPT theorem in two-dimensional
theories of local observables}", Comm. Math Phys. 143 (1992), 315-332
\endREF
\REF(Buch1)Buchholz D., ``{\it On the structure of local quantum
fields with non-trivial interaction}", in: Proc. of the Int. Conf. on
Operator Algebras,  Ideals and their Applications in Theoretical
Physics, Baumg\"artel, Lassner,  Pietsch, Uhlmann, (eds.) pp. 146-153.
Leipzig: Teubner Verlagsgesellschaft 1978.
\endREF
\REF(BuDa1)Buchholz D., D'Antoni C., private communication.
\endREF
\REF(BucSM1)Buchholz D., Schulz-Mirbach H. , ``{\it Haag duality in
conformal quantum field theory}", Rev. Math. Phys. 2 (1990) 105.
\endREF
\REF(Dira1)Dirac, P.A.M., ``{\it Wave equations in conformal space}",
Ann. Math. 37 (1936) 429.
\endREF
\REF(DoLo1)Doplicher S., Longo R. , ``{\it Standard and split
inclusions of von Neumann algebras}", 	Invent. Math. 73 (1984),
493.
\endREF
\REF(Haag1)Haag R. , ``{\it Local Quantum Physics}'', Springer Verlag,
 Berlin Heidelberg 1992.
\endREF
\REF(Fred1)Fredenhagen K. , ``{\it Generalization of the theory of
superselection sectors}", in The algebraic Theory of Superselection
Sectors, D. Kastler ed., World Scientific, Singapore 1990.
\endREF
\REF(FrJo1)Fredenhagen, K., J\"orss M., in preparation, see: J\"orss
M.,  {\it On the
existence of point-like localized fields in conformally invariant
quantum physics}", Desy Preprint.
\endREF
\REF(FrGa1)Gabbiani F., Fr\"ohlich J., ``{\it Operator algebras and
Conformal Field Theory}", pre- print EHT Z\"urich, 1992.
\endREF
\REF(StCe1)Gilman L., Jerison M. , ``{\it Rings of continuous
functions}"
 Springer-Verlag, New York Heidelberg Berlin 1976.
\endREF
\REF(GuLo1)Guido D., Longo R., ``{\it Relativistic invariance and
charge conjugation in quantum field theory}", Comm. Math. Phys. 148
(1992), 521-551. \endREF
\REF(HiLo1)Hislop P., Longo R., ``{\it Modular structure of the local
observables 		associated with the free massless scalar field
theory}", 		Comm. Math. Phys. 84 (1982), 84.
\endREF
\REF(Long1)Longo R., ``{\it Algebraic and modular structure of von
Neumann  algebra of physics}", Proc. of Symposia in Pure Mathematics
38 (1982), part 2, 551-566.
\endREF
\REF(LuMa1)Mack G., L\"uscher M., ``{\it Global conformal invariance
in quantum field theory}", Comm. Math. Phys. 41 (1975) 203.
\endREF
\REF(Post1)Postnikov, ``{\it Le\c cons de geometrie. Groupes et
algebres de Lie}, \'Editios Mir, Moscou, 1985.
\endREF
\REF(Sega1)Segal I. E., ``{\it Causally oriented manifolds and
groups}", Bull. Amer. Math. Soc. 77 (1971) 958.
\endREF
\REF(StraZs1)Str\u atil\u a S., Zsido L. , {\it Lectures on von
Neumann algebras}, Abacus press, England 1979.
\endREF
\REF(StWi1)Streater R.F., Wightman A.S. , {\it PCT, spin and
statistics, and all that}, Addison Wesley, Reading, MA 1989
\endREF
\REF(ToMP1)Todorov I. T., Mintchev M. C., Petkova V. B.,
{\it Conformal invariance in quantum field theory}, Publ. Scuola Normale
Superiore, Pisa, 1978.
\endREF
\REF(Weyl1)Weyl H., {\it Space-Time-Matter}, Dover Pubblications Inc.,
1950.
\endREF
}

\inputreferences

\section Introduction
\par

Haag duality in Quantum Field Theory is the property
that local observable algebras maximally obeys the
causality principle: if $\R(\O)$ is the von Neumann
algebra of the observables localized in the double cone
$\O$ of the Minkowski space $M$, then $\R(\O)$ is the
commutant of the von Neumann algebra $\R(\O')$ of the
observables localized in the space-like complement $\O'$
of $\O$ $$\R(\O')= \R(\O)'$$

Duality plays an important role in the structural
analysis of algebraic Quantum Field Theory [\ref(Haag1)] and
has long been verified in free field models [\ref(Arak1)].
 If the local algebras are generated by Wightman
fields, Bisognano and Wichmann [\ref(BiWi1)] have shown the
general result that duality holds for the von Neumann
algebras associated with  wedge shaped regions $W$,
namely
 $$\R(W')= \R(W)'$$
 where $W$ is any Poincar\'e transformed of the region
$\{ x\in M \mid x_1>|x_0|\}$.
 This property is called essential duality since it allows
to enlarge the original observable algebras of double
cones so that Haag duality holds true.  Their basic
result is obtained by the computation of the
Tomita-Takesaki modular operator $\Delta_W$ associated
to $\R(W)$ with respect to the vacuum vector $\Omega$
[\ref(StraZs1)], the latter being cyclic and separating
because of the Reeh-Schlieder theorem. In this case the
modular group is the (rescaled) one-parameter group of pure
Lorentz transformations leaving $W$ invariant and the
modular conjugation $J_W$ is the product of the PCT symmetry
and a rotation; the essential duality then follows at
once by Tomita's commutation theorem
$$\R(W)'=J_W\R(W)J_W=\R(W').$$ This identification of
the modular group has several interesting consequences;
beside duality, we mention here, as an aside,
 the relation with the
Hawking effect, see [\ref(Haag1)], due to the KMS
(temperature) condition characterizing this evolution, and the
Poincar\'e covariance of the superselection sectors with finite
statistics in [\ref(GuLo1)].
 \par
 One should not expect a sharp geometrical
description for the modular group of the algebra
$\R(\O)$ of a double cone $\O$, since in general there
are not enough space-time symmetries that preserve
$\O$; however, as a consequence of the
Bisognano-Wichmann theorem,  the modular group of
$\R(\O)$ has a geometrical meaning in a conformally
invariant theory [\ref(HiLo1)].

 The purpose of the present work is to provide an
intrinsic, Wightman field independent, algebraic derivation of
the Bisognano-Wichmann theorem in the case of a
conformally invariant theory.

 We were motivated by a recent general result of
Borchers [\ref(Borc1)] showing that part of the geometric
behavior of $\Delta_W$ follows automatically from
the positivity of the energy-momentum operators:
$\Delta_W$ has the expected commutation relations with
the translation operators. This
 result  however does not furnish the commutation
relations between modular operators associated with
different wedges and, in space-time dimension greater
than two, does not provide the Bisognano-Wichmann
identification of the modular automorphism group of
$\R(W)$, indeed simple counter-examples illustrate how it may be
violated in general.

We shall show that essential duality
holds automatically in a conformal theory of any
space-time dimension.
 But Haag duality for double cones fails in general
because the local algebras actually live in a
superworld $\tilde M$ [\ref(Sega1),\ref(LuMa1)], a
$\infty$-sheeted cover of a compactified Minkowski space.
However there is a natural procedure to extend the original net of
local algebras to a causal pre-cosheaf  of
von Neumann algebras (i.e. an inclusion preserving map
$\O\to\R(\O)$) on $\tilde M$, and Haag duality holds there.

As a further consequence, we shall show that an algebraic
conformal theory admits automatically a PCT symmetry [\ref(StWi1)].

In particular for a M\"obius covariant  pre-cosheaf of
local algebras on $S^1$, duality holds on $S^1$ without
any assumption other than positivity of the energy.
 Nevertheless duality fails in general on the cut circle
$\simeq\Re$. This phenomenon, discussed in [\ref(BucSM1)]
in this specialization, is
already present in [\ref(HiLo1)] concerning the time-like duality;
 we shall reformulate the
examples in [\ref(HiLo1)] to get models of conformal  theories on
$S^1$ with the  desired properties.

Our results parallel an independent work of Fredenhagen
and J\"orss [\ref(FrJo1)] where similar results are obtained by a
different route: they construct Wightman fields
associated with an algebraic conformal field theory on $S^1$,
under a finite multiplicity condition for the M\"obius
representation.

Part of our analysis extends to the general case of Poincar\'e
covariant nets. Assuming that the modular group of the von Neumann algebra
of any given wedge region has a geometrical meaning, we can show that the
Bisognano-Wichmann interpretation holds, provided the net fulfills the
split property.
The latter condition is indeed necessary and guarantees the uniqueness
of the Poincar\'e covariant action.

After our work was completed we received a preprint of Gabbiani and
Fr\"ohlich [\ref(FrGa1)] that contains similar analysis for
algebraic conformal field theories on $S^1$.

\section 1. Algebraic Conformal Quantum
Field Theory, general setting.
 \par
 In this section we describe a Conformal Quantum
Field Theory in the algebraic approach.
Our aim is to give a self-contained introduction in our setting
of known features.
\medskip
 \subs{Geometrical preliminaries.}
 In the following we consider a Lie group $G$ acting by local
diffeomorphisms on a manifold $M$, i.e.
there exists an open set $W\subset G\times M$ and a
$C^\infty$ map
 $$\matrix{T:&W&\to&M\cr&(g,x)&\mapsto&T_gx\cr}\eqno(1.1)$$
 with the following properties:
 \item{$(i)$} $\forall x\in M$, $V_x\equiv \{g\in G:(g,x)\in
W\}$ is an open connected neighborhood of the identity $e\in G$
 \item{$(ii)$} $T_ex=x$, $\forall x\in M$
 \item{$(iii)$} If $(g,x)\in W$, then
$V_{T_gx}=V_xg^{-1}$ and moreover for any $h\in G$ such
that $hg\in V_x$
$$T_h T_gx=T_{hg}x$$
 \par
  We say that a local action of a Lie group $G$  on a
 manifold M is {\it quasi-global} if the open set
 $$\{x\in M:(g,x)\in W\}$$
 is the complement of a meager set $S_g$, and the following
equation holds:
 $$\lim_{x\to x_0}T_gx=+\infty,
\qquad g\in G,\quad x_0\in S_g\eqno(1.2)$$
 where $x$ approaches $x_0$ out of $S_g$ and a point goes
to infinity when it is eventually out of any compact
subset of $M$.
\claim{1.1 Proposition} If $T$ is a transitive quasi-global action of $G$ on
$M$, then there exists a unique
``$T$-completion"  of $M$, i.e. a manifold $\overline M$
such that $M$ is a dense open subset of $\overline M$ and
the action $T$ extends to a transitive global action on
$\overline M$.
\par
\proof Let $\E$ be the space of the continuous bounded functions
$\f$ on $M$ such that, for any given $g\in G$,
the function $T^*_g\f$ defined by
 $$T^*_g\f(x)\equiv\f(T_gx),\qquad x\not\in S_g$$
 has a continuous extension to $M$.
 \np
 Since a countable union of meager sets is meager,
it is easy to see that $\E$ is a $C^*$-algebra, and $g\to T^*_g$ is an
action of $G$ by isomorphisms of $\E$.
 \np
 Moreover $\E$ contains the algebra $C_0(M)$ of
continuous functions on $M$ vanishing at infinity because, by
condition (1.2),
 $$T^*_g\f(x)\equiv\cases{\f(T_gx)&if $x\not\in S_g$\cr 0&if
$x\in S_g$\cr}$$
 is a continuous function, for any $\f\in C_0(M)$ and for any fixed
$g\in G$. We denote by $C_G$ the minimal $C^*$-subalgebra of $\E$ containing
$C_0(M)$ and globally invariant under the action $T^*$.
 \np
 We show that the spectrum $\overline M$ of $C_G$ is the requested
completion, where the (global) action of $ G$ on $\overline
M$ is given by
 $$T_gp=p\cdot T^*_g,\qquad p\in\overline M.$$
 We notice that $C_G$ does not necessarily contain the
identity, therefore $\overline M$ is not compact in general.
 Moreover, by well known arguments [\ref(StCe1)], the spectrum of the
$C^*$-algebra obtained adding the identity to $C_G$ is a
compactification of $M$ and is the one-point compactification of $\overline M$,
hence the natural embedding $M\hookrightarrow\overline M$ is dense.
 \np
 The transitivity of the action on $\overline M$ follows by the minimal choice
of $C_G$. In fact, let us consider the  orbit in $\overline M$ containing $M$,
 $$\overline M_0\equiv\{T_gx:g\in G\}$$
 where $x$ is a point of $M$ as a subspace of $\overline
M$. By the transitivity on $M$, $\overline M_0$ does not
depend on $x$ and the action of $ G$ is transitive on
$\overline M_0$. We have natural embeddings
 $$ C_0(M)\subseteq C_0(\overline M_0)\subseteq C_G$$
and $C_0(\overline M_0)$ is globally invariant under $T^*$. Since
$C_G$ is, by construction, the minimal
$C^*$-algebra with this properties, then
$\overline M_0=\overline M$.
 \np
 By minimality, any other $T$-completion of $M$ contains $\overline M$, and
therefore the transitivity requirement implies uniqueness.
 \np
 Finally, $\overline M$ is a homogeneous space for the Lie group $G$, hence it
is a $C^\infty$ manifold.
 \par
 \claim{1.2 Proposition} In the hypotheses of the
previous proposition, we consider the universal
covering $\tilde G$ of $G$ and the universal covering $\tilde M$
of $\overline M$. Then the action $T$ lifts to a
transitive global action $\tilde T$ of $\tilde G$ on $\tilde
M$.
 \par
 \proof The group $\tilde G$ has a canonical global
transitive action on $\overline M$, i.e. $T\cdot\p$ where
$\p:\tilde G\to G$ is the covering map. The identity
component of the isotropy group of a point $(\tilde
G_x)_0$ is a connected closed subgroup of $\tilde G$,
therefore [\ref(Post1)], the manifold
 $$\tilde M\equiv\tilde G/(\tilde G_x)_0$$
 is  simply connected and is a covering of $\overline M$,
i.e. it is the universal covering of $\overline M$. Then $\tilde M$ is a
homogeneous space for $\tilde G$, and the thesis follows.
 \par
 \subs{Conformal action on the Minkowski space.} In the
following we specialize the preceding analysis to the action
of the conformal group on the Minkowski space $M$.
The action of the conformal group $\C$ on a point of $M$
(resp. of the universal covering $\tilde\C$ on $\tilde M$) will be denoted
by $x\to gx$.
 As is known, $M$ is the manifold $\Re^d$ with constant
pseudo-metric tensor $q$ with signature
 $$\matrix{+1,&-1,\dots,-1\cr&$\upbracefill$\cr
 &{\rm(d-1)-times}\cr}\qquad.\eqno(1.3)$$
The corresponding Minkowski norm of a vector $x\in M$ is
$$x^2=x_0^2-x_1^2-\dots-x_{d-1}^2.$$
When $d>2$ the conformal group consists of the local
diffeomorphisms $\f$ of $M$ which preserve the
pseudo-metric tensor $q$ up to a non-vanishing
function $\m$:
$$\f^*q=\m q$$
 This group is a ${(d+2)(d+1)\over 2}$ dimensional Lie  group, and
the Lie algebra of its identity component $\C$ is generated
(as a vector space) by the following  objects:
\medskip
\settabs 2\columns
\+\qquad translations&$d$ generators\hfill$(1.4a)$&\cr
\+\qquad boosts&$d-1$ generators\hfill$(1.4b)$&\cr
\+\qquad rotations
 &${(d-1)(d-2)\over2}$ generators\hfill$(1.4c)$&\cr
\+\qquad dilations&1 generator\hfill$(1.4d)$&\cr
\+\qquad special transformations
 &$d$ generators\hfill$(1.4e)$&\cr
 \medskip\noindent
The special transformations are the elements of the form
$$\r\t_a\r$$
where $\t_a$, $a\in \Re^d$, is a translation and $\r$ is the
relativistic ray inversion: $$\r x=-{x\over x^2}\qquad x^2\not=0.$$
  \par When $d\leq2$ the conformal group $\C$ is,
by definition, the one generated by the
transformations described in (1.4).
 \par
 \subs{ The Conformal Universe.}
 The conformal group $\C$ acts quasi-globally on $M$, and
therefore the completion $\overline M$ is defined.
 \par
 For convenience of the reader, we follow [\ref(ToMP1)] and
give the more explicit construction of $\overline M$ by Dirac
and Weyl [\ref(Dira1),\ref(Weyl1)]. Given $\Re^{d+2}$
with signature
$$\matrix{+1,&-1,\dots,-1&+1\cr&$\upbracefill$&\cr
 &d-{\rm times}&\cr}$$
 we consider the
manifold $N$ whose points are the isotropic  rays of the
light cone, i.e.
$$N=\{(\x_0,\dots \x_{n+1})\in\Re^{d+2}\backslash\{0\}:
+\x_0^2-\x_1^2-\dots -\x_d^2+\x_{d+1}^2=0\}/\Re^*$$
where $\Re^*=\Re\backslash\{0\}$ acts by multiplication on $\Re^{d+2}$.
The Lie group $PSO(d,2)$ acts transitively on this
manifold by global diffeomorphisms, and it is easy to
check that the map $M\to N$ given by
$$
\left\{\eqalign{\x_i&=x_i\quad i<d\cr
       \x_d&={1-x^2\over 2}\cr
       \x_{d+1}&={1+x^2\over 2}\cr}\right.\eqno(1.5)
$$
is a dense embedding such that the restriction of
the action of $PSO(d,2)$ to $M$ corresponds to the
conformal transformations. The uniqueness proven in
Proposition 1.2 shows that $N=\overline M$.
\par
The Dirac-Weyl description shows  $\overline M$ to be
 a compact manifold diffeomorphic to
$(S^{d-1}\times S^1)/\Ze_2$.
When $d>2$, the universal covering $\tilde M$ is infinite
sheeted and is diffeomorphic to $(S^{d-1}\times\Re)$.
The map described in equation (1.5) lifts to a
natural  embedding of $M$ into $\tilde M$, and the covering
map from $\tilde M$ to $\overline M$ shows that $\tilde M$
contains infinitely many copies of $M$ as submanifolds.
\par
One of the main advantages of dealing with the manifold
$\tilde M$ instead of $\overline M$ is that a global
causal structure is naturally defined on it, i.e. the
time ordering gives rise to a global ordering relation
which extends the ordering on $M$, and a notion of
(non-positive definite) geodesic distance is also
well-defined, and is locally equivalent to the one in
$M$. As a consequence each point divides $\tilde M$ in
three parts, the relative future, i.e. the
points at time-like distance which follow the point, the
relative past, i.e. the points at time-like distance
which precede the point, and the relative present, i.e.
the points at space-like (or light-like) distance from
the point. As we shall see, conformal quantum field theories live naturally on
$\tilde M$.
 \par
 When $d\leq2$, $\overline M$ is diffeomorphic to $(S^1)^d$, and therefore its
 universal covering is $\Re^d$, but this manifold is rather unphysical.
 In fact we may find two space-like separated embeddings of $M$ in $\Re^d$.
 Therefore we use the convention $\tilde M\equiv S^1\times\Re$ when
$d=2$, and $\tilde M\equiv S^1$ when $d=1$. \par
We also mention that when $d$ is odd ($d\not=1$) the manifold
$\overline M$ is not orientable, and physical theories
live on orientable coverings of $\overline M$ (cf.
Proposition 1.7 and the odd-dimensional examples at the end of Section 2).
 \par
 In the following we shall consider the
family $\tilde\K$ of the subregions of $\tilde M$ which are
images of  double cones in $M$ under  conformal transformations in
$\tilde\C$.  We notice that all double
cones, wedges and light-cones of $M$  belong to the family
$\tilde\K$, see e.g. [\ref(HiLo1)].
Now we list some properties of $\tilde\K$:
 \claim{1.3 Proposition}
 \item{$(i)$} All elements of $\tilde\K$ are open
contractible precompact submanifolds of $\tilde M$. They are a fundamental set
of neighborhoods for $\tilde M$.
 \item{$(ii)$} $\tilde\C$ acts ``transitively" on
$\tilde\K$, i.e.
 $$\forall\ \O_1,\O_2\in\tilde\K\quad\exists
g\in\tilde\C:g\O_1=\O_2$$
 \item{$(iii)$} The identity component $\tilde\C(\O)_0$
of the group of the conformal diffeomorphisms that preserve $\O$,
 $$\tilde\C(\O)=\{g\in\tilde\C:g\O=\O\}$$
 acts transitively on $\O$, $\O\in\tilde\K$.
 \item{$(iv)$} The space-like complement $\O'$ of a region
$\O\in\tilde\K$ belong to $\tilde\K$.
 \item{$(v)$} The family $\tilde\K$ is not a net, in fact the union of a region
and of its causal complement is not contained in any region of $\tilde\K$.
 \par
\proof Immediate.
 \par
 With each region $\O\in\tilde\K$, we shall associate a one-parameter group
$\L^\O_t$ of conformal
transformations which preserve $\O$ and commute with all $\O$-preserving
conformal transformations:
 $$\O\in\tilde\K\to\{\L^\O_t,\quad t\in\Re\}\subset\tilde\C$$
 These groups will have the
following coherence property:
 $$\L^{\O_2}_t=g^{-1}\L^{\O_1}_tg,\qquad\O_2=g\O_1,\,
\O_i\in\K\quad g\in\tilde\C$$
 Therefore, by Proposition 1.3$(ii)$, they will be completely determined
if we assign $\L^\O_t$ for one region $\O$.
 For reader's convenience, we describe explicitly the conformal
transformations for three particular regions in $M$ (cf.
[\ref(HiLo1),\ref(Buch1)]).
 \np
 The wedge $W_1$:
 $$W_1=\{x\in M:x_1>|x_0|\}.$$
 The group $\L^{W_1}_t$ is the one-parameter group of pure Lorentz
transformations (boosts) along the $x_1$ axis; its action on $(x_0,x_1)$ is
given by the matrices
 $$\left(\matrix{\cosh 2\p t&-\sinh2\p t\cr -\sinh2\p t
&\cosh2\p t\cr}\right)$$
 The double cone $\O_1$:
 $$\O_1=\{x\in M:|x_0|+|\vec x|<1\}$$
 The group $\L^{\O_1}_t$ commutes with the rotations, hence is determined by
its action on the $(x_0,x_1)$-plane:
 $$\L^{\O_1}_tx_{\pm}={(1+x_{\pm})-e^{- 2\pi t}(1-x_{\pm})\over
   (1+x_{\pm})-e^{- 2\pi t}(1+x_{\pm})}$$
 where we posed $x_{\pm}=x_0 \pm x_1$.
\np
 The future cone $V_+$:
 $$V_+=\{x\in M:x_0>0,\,x^2>0\}$$
 The one-parameter group $\L^{V_+}_t$ is the dilation subgroup:
 $$\L^{V_+}_t=D(e^t)$$
 \medskip
 \subs{ The universal covering $\tilde\C$.}
The universal covering $\tilde\C$ of the conformal group
$\C$ turns out to be a central extension of $\C$ with fiber
$\Ze\times\Ze_2$. The $\Ze_2$ component acts trivially on
$\tilde M$ and $\tilde\C/\Ze_2$ acts effectively on it
(see e.g. [\ref(ToMP1)]).
 \par
 Now we prove some simple properties on the conformal
group we shall need in Section 2.
\claim{1.4 Proposition} The groups $\C$ and $\tilde\C$ are
perfect groups, i.e. they coincide with their commutator subgroups.
 \par
 \proof Since both $\C$ and $\tilde\C$ are semi-simple Lie
groups, the result follows by the observations in [\ref(Post1)], p. 345.
 \par
\claim{1.5 Proposition} The conformal transformations $R_i$,
$i=1,\dots,d-1$ given by
$$R_ix=-{1\over x^2}(x_0,\dots
x_{i-1},-x_i,x_{i+1},\dots,x_{d-1})$$
are in the identity component of $\C$.
Each $R_i$ has only two liftings $\pm R_i$ of order 4 in
$\tilde\C$.
 \par
 \proof (Sketch) We observe that
$$U(\a)\equiv\t_i(-\cot\a)D((\sin\a)^{-2})R_i\t_i(-\cot\a)$$
is a one-parameter subgroup of $\C$, where $\t_i(\cdot)$ are the
translations along the $i-th$ axis. Moreover $U(\a)$ satisfies
$$\left\{\eqalign{U(0)&=U(\p)=e\cr
U({\p\over2})&=R_i\cr}\right.,$$
hence the first assertion follows.
Lifting this subgroup to $\tilde U(\a)\in\tilde\C$ we get a
group of period $2\p$, therefore $\pm R_i\equiv\tilde
U(\pm{\p\over2})$ are the requested liftings. Since the fiber
of the covering $\tilde\C\to\C$ is $\Ze\times\Ze_2$, any
other lifting has infinite order.
 \par
\claim{1.6 Proposition} For each $i=1,\dots,d-1$, the translation
subgroup and  $R_i$ generate the conformal group $\C$. The
same result holds in $\tilde\C$ when $R_i$ is
replaced by any of its liftings $\pm R_i$
(cf. [\ref(ToMP1)] for an analogous statement in the
4-dimensional case).
 \par
 \proof (Sketch) A straightforward calculation shows that the
equation
 $$\t_i(a)R_i\t_i(1/a)R_i\t_i(a)R_i=D(a^2)$$
holds in $\C$.
Then, since dilations are the transformations $\L^{V_+}$ for the
future light cone and such a cone can be mapped into any wedge
by a suitable product of
 translations and $R_i$, boosts are in the group generated by $\t$ and $R_i$
[\ref(HiLo1)]. The first statement follows because boosts
and translations generate the Poincar\'e group and, together with
$R_i$, the conformal group $\C$ (cf. formulas (1.4)).
The statement for  $\tilde\C$ is proven in a similar way.
 \par
 \claim{1.7 Proposition} When $d$ is odd, the change of
sign of a space coordinate $P_i$ is in the identity
component of $\C$. As a consequence $\overline M$ is not
orientable.
 \par
 \proof The first part of the Proposition follows from
straightforward calculations similar to those in
Propositions 1.5 and 1.6. The rest follows because the Jacobian of $P_i$ is
negative.
\par
 \medskip
\subs{ Conformal Quantum Field Theories.}
A local Conformal Quantum Field Theory in dimension $d$
is described by a causal additive pre-cosheaf of
von Neumann algebras on the double cones of $M=\Re^d$ with
Minkowski structure, i.e. a map
$$\A:\O\to\A(O),\qquad \O\in\K$$
where $\K$ is the family of the double cones in $M$, such that
$$\eqalign{ \O_1\subset\O_2&\imply\A(\O_1)\subset\A(\O_2)\cr
\A(\O)&\subset\A(\O')'\hskip3.cm {\rm(causality)}\cr
\A(\cup_n\O_n)&=\vee_n\A(\O_n)\hskip2.5cm {\rm(additivity)}}$$
where $\O'$ is the space-like complement of $\O$ and
the regions $\O$, $\O_n$ and $\cup\O_n$ belong to $\K$.
The $\A(\O)$ are supposed to act on a common Hilbert space
$\H$. Since the family $\K$ is a direct set, the map $\O\to\A(\O)$ is indeed a
net and the quasilocal $C^*$-algebra $\A_0$ is defined as
the direct limit of the local algebras. The algebras associated with
general open regions in $M$ are defined  by additivity.
 \par
 Now we describe the conformal covariance assumption. As already explained, the
group $\C$ acts locally on $M$.
 We observe that, since double cones are precompact, using the
notation of formula (1.1) there exists an open neighborhood $V_\O$ of
of the identity in $\tilde\C$ such that $V_\O\times\O\subset W$,
 and therefore the elements in $V_\O$ give rise to diffeomorphisms of $\O$ into
$M$.
 \par
 We assume that $\C$ acts  locally by covariant
automorphisms of the pre-cosheaf $\A$, namely, for any
$\O\in\K$ there is a weakly continuous map from the open set
$V_\O$ to the set $Iso(\A(\O),\A_0)$ of isomorphisms of
$\A(\O)$ into $\A_0$,
 $$g\to\a_g^\O$$
  with the following properties.
\np
If $\O_1\subset\O_2$ are double cones,
 $$\a_g^{\O_2}\big|_{\A(\O_1)}=\a_g^{\O_1}$$
 Because of the preceding property, and since double cones form
a net, we may drop the superscript which specify the region,
$\a_g\equiv\a_g^\O$.
\np
 The map $g\to\a_g$ is a local action, i.e.
 $$\a_{hg}=\a_h\cdot\a_g$$
 when it makes sense.
\np
 The local action is covariant, i.e.
$$\a_g\A(\O)=\A(g\O).$$ \par
Finally we assume the existence of a local unitary representation $U$ of $\C$
and of a $U$-invariant vector $\Q$, cyclic for $\cup\A(\O)$,  such that
$$U(g)A\Q\equiv \a_g(A)\Q\qquad A\in\A(\O),\quad g\in
V_\O\eqno(1.6)$$
The generators of the (local)
one-parameter subgroups are well-defined selfadjoint operators
on $\H$.
The energy-momentum is assumed to be positive.
 \par
 While the usual energy $H$ corresopnds to the Lie algebra
generator $h$ of the time translations, the conformal energy $K$
corresponds to the Lie algebra element
 $$k\equiv h+\r h\r\eqno(1.7)$$
 where $\r h\r$ is the adjoint action of $\rho$ on $h$.
 It is clear that if $H$ is positive, $K$ is the sum of positive
operators, and therefore is positive too. The converse is also
true, and can be checked on explicit realizations of the irreducible
positive-energy unitary representation of $\C$ [\ref(Sega1)].
 \par
 By the Reeh-Schlieder theorem $\Q$ is cyclic and separating for the algebras
$\A(\O)$ and $U$ can be defined directly from (1.6).
Now we show how the pre-cosheaf $\A$ may be canonically
extended to the superworld $\tilde M$.
 \par
 We observe that the map $g\to U(g)$ may be considered as a
 local unitary representation of the universal covering
$\tilde\C$, where  $g\in\tilde V_\O$, the identity component of
the pre-image of $V_\O$ under the covering map. Then, since
$\tilde\C$ is simply connected, $U$ extends to a global unitary
representations of $\tilde\C$ on $\H$.
 \par
 \claim{1.9 Lemma} If $\O\in\K$, $g\in\tilde\C$ and
$g\O=\O$ then $U(g)\A(\O)U(g)^*=\A(\O)$.
\par
\proof We take $x\in\O$, then $gx\in\O$ and since the
connected component of the stabilizer of $\O$ acts
transitively on $\O$ (see Proposition 1.3) we find
$h\in\tilde\C(\O)_0$ such that $h^{-1}gx=x$.
 Then $h^{-1}g\in\tilde\C_x$, which is connected by construction
(cf. the proof of Proposition 1.2).
 Therefore we may find a one-parameter family
$l(t)\in{\tilde\C}_x$ such that $l(0)=id$ and
$\l(1)=h^{-1}g$, and there exists a neighborhood
$B_x\in\tilde\K$ of $x$ such that $l(t)B_x\subset\O$,
$t\in[0,1]$, hence the local covariance of $\a$ implies covariance for all
$t\in[0,1]$, i.e.
 $$\a_{h^{-1}g}\A(B_x)=\A(h^{-1}gB_x).$$
 Since $h\in\tilde\C(\O)_0$, we may connect it with $id$ staying
inside $\tilde\C(\O)_0$, and  the same argument used
before implies $\a_h\A(\O)=\A(\O)$. As a consequence,
$$\a_g\A(B_x)=\a_h\a_{h^{-1}g}\A(B_x)=\A(gB_x)\subset\A(\O)$$
By additivity $\A(\O)=\vee_{x\in\O}\A(B_x)$, therefore
$$\a_g(\A(\O))=\bigvee_{x\in\O}(\A(gB_x))=\A(g\O).$$
\par
 Since $\tilde\C$ acts transitively on
$\tilde\K$, we may define
  $$\tildA(g\O)=U(g)\A(\O)U(g)^*,\quad g\in\tilde G,\
\O\in\K.\eqno(1.7)$$
 By Lemma 1.9, the map $\tildA$ is well defined on
$\tilde \K$, and defines a pre-cosheaf of
$C^*$-algebras on $\tilde\K$. By definition $\tilde\C$
acts globally covariantly on the pre-cosheaf $\tildA$.
 We observe that if $W\subset M$ is a wedge, then $\tildA(W)$ turns out to be
weakly closed, while $\A(W)$ is not. However $\tildA(W)=\A(W)''$,
therefore, since $\A(W)$ is defined by additivity and $\A$ is causal,
$\tildA(W)\subset\tildA(W')'$. Then, by covariance and transitivity on
$\tilde\K$, causality holds for $\tildA$. The additivity property for
$\tildA$ follows by the analogous property of $\A$.
 Thus we have proven the following:
\claim{1.10 Proposition} The pre-cosheaf $\A$
 extends to a unique, causal, additive pre-cosheaf $\tildA$
on $\tilde M$ and the local action $\a$ canonically
extends to a globally covariant action of the group
$\tilde\C$
 \par
 We have therefore shown the equivalence between the locally covariant
picture, on the manifold $M$, and the globally covariant picture, on
$\tilde M$, for  conformally covariant field theories.
 \par
  As it is explained in ([\ref(Fred1)], see also [\ref(GuLo1)]) we
may define a universal $C^*$-algebra $\tildAnul$, and the
isomorphisms $\ad U(g)$ extend to automorphisms
$\tilde\a_g$ of $\tildAnul$.
 The $C^*$-algebra $\tildAnul$ is larger than $\A_0$ since it contains the
von~Neumann algebras associated with wedge regions, and it is not
necessarily faithfully represented in the vacuum representation [\ref(Fred1)].
 \par
 Finally we mention that explicit models may live on finite coverings of
$\overline M$ or even on $\overline M$ itself (see examples in Section 2).
 In the last part of the
following section we shall illustrate with some examples this phenomenon.
 \section 2. Duality property and the Bisognano-Wichman theorem for conformal
theories.
 \par
In this Section we prove that essential duality holds for an algebraic
conformal field theory on the Minkowski space $M$, and duality for double cones
(and conformally equivalent regions) holds for the
corresponding pre-cosheaf extension on $\tilde M$. Moreover  the
modular unitary group of a region $\O\in\tilde\K$ coincides with a
one-parameter subgroup of the conformal group, namely we have an algebraic
derivation of Bisognano-Wichmann theorem in this case.
 \par
As was explained in the preceding section, the Minkowski space $M$
is embedded in a canonical way in the superworld $\tilde M$, and a
(locally)  conformally covariant pre-cosheaf of von~Neumann
algebras  on $M$ extends uniquely to a pre-cosheaf on $\tilde M$
which is globally covariant with respect to the universal covering
$\tilde\C$ of the conformal group $\C$. Therefore, in the
following, $M$ will always be thought  as a submanifold of $\tilde
M$, and $\A$ as a sub-pre-cosheaf of $\tildA$. We shall denote by
$\R(\O)$ the von~Neumann algebra associated with $\O$, both when $\O\in\K$ and
when $\O\in\tilde\K$
$$\R(\O)=\tilde\A(\O).$$
 \par
 In two-dimensional space-time theories we make the further
assumption that parity is implemented, i.e. there exists a
selfadjoint unitary $U(P)$ such that
 $$U(P)\R(\O)U(P)=\R(P\O)$$
 where $P$ is the change of sign of the space coordinate in
$M$, or the corresponding transformation in $\tilde M$.
 In particular, for chiral theories, namely pre-cosheaves which split into a
product of two pre-cosheaves on $S^1$, parity correspond to the flip
automorphism of the tensor product.
 \par
 We recall that an {\it internal symmetry} of the theory is an automorphism of
the pre-cosheaf, i.e. a consistent family of automorphisms
 $$\g^\O\in Aut(\R(\O)),\qquad \O\in\tilde\K$$
 \claim{2.1 Lemma} Let $\O$ be a region
in $\tilde\K$, $\D_\O$ the modular operator for $(\R(\O),\Q)$
and $U_\O(t)$ the unitaries associated with the transformations
$\L^\O_t$ defined in (1.6). Then $$z(t)=\D_\O^{it}U_\O(-t)$$
is a one-parameter unitary group that commutes with the unitary
representation of $\tilde\C$, implements internal symmetries,
and does not depend on the region $\O$.
 \par
 \proof $d>2.$\quad Let us consider the wedge $W_1$ and the
equation
 $$\D_{W_1}^{it}U(g)\D_{W_1}^{-it}=
U_{W_1}(t)U(g)U_{W_1}(-t).\eqno(2.1)$$
 If $g$ is a conformal diffeomorphism which preserves $W_1$, then
$g$ commutes with $\L^{W_1}_t$. Moreover
$U(g)$ implements an automorphism of $\R(W_1)$, hence, by
Tomita-Takesaki theory (see e.g. [\ref(StraZs1)]),  it commute with
$\D_{W_1}^{it}$ and (2.1) holds. This is also  the case when $g$ is a
translation
along the axes $x_2,\dots,x_{d-1}$ or the trasformations $\pm R_2$
introduced in Lemma 1.5. By Borchers theorem [\ref(Borc1)], $(2.1)$ holds
if $g$ is a translation along $x_0$ or $x_1$, and therefore, by Lemma 1.6,
for all elements of the conformal group.
 As a consequence, $z(t)=\D_{W_1}^{it}U(-t)$
commutes with $U(\tilde\C)$ and therefore, since it preserves $\R(W_1)$, it
preserves $\R(\O)$ for any region $\O\in\tilde\K$.  Hence $z(t)$
is an internal symmetry and commutes with the modular groups of all regions in
$\tilde\K$.  The independence from the region follows by the mentioned
commutation relations. \np
 If $d=2$,   the proof goes on as before provided we substitute
$R_2$ with $R_1P$.
\np
 If $d=1$ (i.e. on $S^1$), we consider the tensor product of
$\tildAnul$ with itself, and we get a chiral  two-dimensional
theory.  Then formula (2.1) holds for this theory, hence for the
original one-dimensional theory.
 \par
 \claim{2.2 Corollary} The given
unitary representation $\{U(g):g\in\tilde\C\}$ is the unique
representation of the conformal group which implements a covariant
action on the pre-cosheaf $\tildA$ and preserves the vacuum vector.
 \par
 \proof By Lemma 2.1 we get the equality of the
multiplicative commutators
 $$[\D_{\O_1}^{it},\D_{\O_2}^{is}]=
[U(\L^{\O_1}_t),U(\L^{\O_2}_s)].$$
 Therefore the representation $U$ is intrinsic for $g\in [\tilde\C,\tilde\C]$,
since it it determined by the modular group.
 Since $\tilde\C$ is perfect (Proposition 1.4), $U$ is completely
determined by the modular operators, hence it is unique.
 \par
 \claim{2.3 Theorem}  Let $\A$ be a conformally covariant
pre-cosheaf on $M$, $\tilde A$ the corresponding pre-cosheaf on
$\tilde M$. Then:
 \item{$(i)$} Essential duality holds for the
pre-cosheaf $\A$, and duality holds for the pre-cosheaf
$\tildA$.
 \item{$(ii)$} If $\O$ is a region in $\tilde\K$, $\D_\O$
the modular operator for $(\R(\O),\Q)$ and $U_\O(t)$
the unitary associated with $\L^\O_t$, then
$$\D_\O^{it}=U_\O(t).$$
\par
\proof $(i)$ Let us consider the wedge $W_1$. By causality,
$\R(W'_1)$ is a subalgebra of $\R(W_1)'$, and is globally
stable under the action of the modular group of $\R(W_1)$ because
the latter acts geometrically by Proposition 2.1.  Therefore, by
Tomita-Takesaki theory and the cyclicity of the vacuum,
it coincides with $\R(W_1)'$. Since the
Poincar\'e group has a global causally-preserving action  on $M$
the relation
 $$\R(W')=\R(W)'$$
 holds for each wedge $W$, i.e. essential duality holds for the pre-cosheaf
$\A$. Since $\tilde\C$ has a global causally-preserving action  on $\tilde
M$ the relation
 $$\R(\O')=\R(\O)'$$
 holds for any region $\O\in\K$, i.e. duality holds for the pre-cosheaf
$\tildA$.
 \np
 $(ii)$ We have to show that $z(t)$ in Lemma 2.1 vanishes.
Let us consider one of the conformal transformations $\pm R_1$
introduced in Proposition 1.5.
By definition, $R_1 W_1=W'_1$. Then, by essential duality,
 $$U(R_1)\D_{W_1}U(R_1)=\D_{W_1}^{-1}.$$
 Moreover,
 $$U(R_1)U_{W_1}(t)U(R_1)=U(R_1\L_t^{W_1}R_1)=U_{W_1}(-t).$$
 Therefore,
 $$z(t)=U(R_1)z(t)U(R_1)=U(R_1)\D_{W_1}^{it}U_{W_1}(-t)U(R_1)=z(-t),$$
hence $z(t)=I$.
 \par
Now we show that in even dimensions there exists a
canonical antiunitary $\T$ which implements the PCT
transformation on $\tilde M$, i.e. implements a
conformal transformation of $\tilde M$ that restricts to
change of sign of all coordinates on $M$:
 $$\T\R(\O)\T=\R(\b\O),\qquad\O\in\K$$
 where $\b x=-x$, $x\in M$.
 Such an antiunitary is not unique, since $\T V$ still
implements a PCT whenever $V$ is a self-adjoint  unitary which
implements a internal symmetry. A suitable
(positivity) condition will fix a canonical choice for $\T$.  We
choose the wedge $W_1$ as in formula (2.1) and consider the
unitary $S_{W_1}$ which corresponds to the change of sign of the
coordinates $x_2,\dots x_{d-1}$ (which is an element in the
identity component of $\tilde\C$ when $d$ is even). Then,
with each wedge $W$ in the family $\W_0$ of the
Lorentz transformed regions of $W_1$, we associate the unitary
$S_W\equiv U(g)S_{W_1}U(g^{-1})$ where $gW_1=W$, $g\in\Lor$.
 \claim{2.5 Theorem} The anti-unitary $\T=J_W S_W$,
$W\in\W_0$, implements a PCT transformation and does not
depend on the choice of the wedge $W\in\W_0$. It is the unique anti-unitary
such that
 $$\T\R(\O)\T=\R(\b\O)\qquad\O\in\K$$
and
 $$(\Q,A\T S_W A\Q)\ge 0\quad\forall
A\in\R(W),\quad W\in\W_0\eqno(2.2)$$
 \par
\proof The same arguments used in the proof of Theorem 2.1
imply that
 $$J_{W_1}U(g)J_{W_1}=U(r_1gr_1),\quad g\in\C$$
where $J_{W_1}$ is the modular conjugation for the algebra
$\R(W_1)$ and $r_1$ is the reflection with respect to the edge of the wedge.
Since $J_{W_1}\R(W_1)J_{W_1}=\R(W'_1)=\R(r_1W_1)$ we get
$$J_{W_1}\R(\O)J_{W_1}=\R(r_1\O),\quad\forall\O\in\K.$$
Hence $J_{W_1} S_{W_1}$
implements the transformation $\b$ defined by $\b x=-x$.
Therefore the map
 $$g\to\T U(\b g\b)\T\qquad g\in\tilde\C$$
 is a covariant unitary representation of $\tilde\C$.
 By the uniqueness proved in Corollary 2.2, and since $\b$
commutes with the Lorentz group $\Lor$ we get
 $$\T U(g)\T=U(g)\qquad g\in\Lor.$$
 \np
The previous commutation relation implies that $\T$ does
not depend on $W$:
\np
in fact given two wedges $W_1$, $W_2$ we can find an
element $g\in\Lor$ such that $gW_1=W_2$,  therefore we
have
 $$J_{W_1}S_{W_1}=U(g)J_{W_2}U(g)U(g)^*S_{W_2}U(g)^*=
J_{W_2}S_{W_2}.$$
 \np
 Since  $\T$ implements the transformation $\b$, $\T S_W$ is
an antilinear conjugation which maps $W$ in $W'$, hence, by
essential duality, $\T S_W\R(W)\T S_W=\R(W)'$.
Therefore condition (2.2) is the characterization of
the  modular conjugation $J_W$ [\ref(StraZs1)], and the uniqueness follows.
 \par
We observe that in odd dimensions we may implement a PCT
transformation up to a change of sign of one space coordinate.
Therefore a complete PCT is implemented if and only if such a
change of sign is unitarily implemented.
\rmclaim{2.6 Remark} The split property [\ref(DoLo1)]
holds automatically for a
conformal pre-cosheaf on $(S^{d-1} \times S^1 ) /\Ze_2$, where it is
equivalent to the distal split property [\ref(DoLo1)]
 by conformal invariance, provided the conformal Hamiltonian
$K$ (that has spectrum equal to $\Na$) has the multiplicity of its eigenvalues
growing at most exponentially. Indeed $e^{-\beta K}$ is a trace-class
operator (for $\beta$ large enough) and nuclearity holds [\ref(BuDa1)].
 \par
\claim{2.7 Corollary} With the assumptions of the previous remark, any
local algebra $\R(\O),\quad \O\in\tilde\K,$ is the unique injective
factor of type $III_1$.\par
\proof Immediate, see [\ref(Long1)].\par
 \subs{Examples derived from free fields.}
We illustrate the present setting by some examples given by the net $\A$ of
local algebras associated with the free massless scalar field on the
$d$-dimensional Minkowski space discussed in [\ref(HiLo1)], with a
reinterpretation of the duality behavior.
 \par
 Recall that the group of conformal time translations (i.e. the
subgroup generated by the conformal energy $K$) has period $2\p$, i.e. the
operator $e^{2\p iK}$ acts as the identity on $\tildA$.
 The $\Ze$-component of the center of $\tilde\C$ acts as an
``helicoidal shift" on $\tilde M$ and that the $2\p\Ze$ time translations
correspond to the even part of the $\Ze$-component of the center of
$\tilde\C$.
Therefore all even copies of $M$ may be identified, and all odd copies
of $M$ may be identified, hence free massless scalar fields lives on the
compactification $\overline M$ or on its 2-fold covering $\overline M_2$. In
the following, if $\O$ is a region in $M$, we shall indicate  with $\O'$ its
space-like complement in $M$, with $\O^t$ its time-like complement in $M$, and
with $\O^c$ its  causal complement in $\overline M$, resp. $\overline M_2$.
\par
\item{$\bullet$}$d$ even.\quad Since the relativistic ray inversion $\r$ is
unitarily implemented, the theory lives on $\overline M$.
 The space-like complement $\O'$ and the time-like complement $\O^t$ of
a double cone $\O\in M$ form a connected region $\O^c$ in $\overline M$, the
causal complement of $\O$ in $\overline M$. The duality property in  $\overline
M$ means
 $$\R(\O)'=\R(\O^c).$$
 Since $\R(\O')$ and $\R(\O^t)$ are subalgebras of $\R(\O^c)$, space-like and
time-like commutativity hold for the pre-cosheaf $\A$.
 Moreover $\R(\O')=\R(\O^c)$, therefore Haag duality holds in $M$.
 If $d>2$, $\R(O^t)\not=\R(\O^c)$, i.e. time-like duality does not
hold in $M$ [\ref(HiLo1)].
 If $d=2$, $\R(\O^t)=\R(\O^c)$, i.e. time-like duality holds in $M$.
 \par
 \item{$\bullet$} $d$ odd ($d\not=1$).\quad Since two different liftings of
 $\rho$ in $\tilde\C$ are unitarily implemented, the theory lives on the
2-fold covering $\overline M_2$.
 The product of the two liftings of the relativistic ray inversion is the
transformation $\O\to\tilde\O$ which maps a region in one of the copies of $M$
in $\overline M_2$ onto the corresponding region in the other copy.
 The causal complement in $\overline M_2$ of a double cone $\O\in M$ is the
region
 $$\O^c=\O'\cup \tilde{\O^t}.$$
 Since
 $$\R(\O')=\R(\O^c)$$
 then Haag duality follows by the duality property for the pre-cosheaf
$\tildA$.
 The ``twisted time-like commutativity" in [\ref(HiLo1)] is the inclusion
 $$\R(\tilde{\O^t})\subset\R(\O)'$$
 i.e. the algebra associated to the twisted time-complement $\tilde{\O^t}$
commutes with $\R(\O)$.
 \par
 \item{$\bullet$} $d=1$.\quad Fix $n>2$ and let $\R(\O)$ be the net on the
$n$-dimensional Minkowski space described above. Then set
 $$\A_n(I)=\R(\O)$$
 where $I\subset\Re$ is an interval of the time axis and $\O$ is the double
cone obtained by causal completion from $I$.
\par
If $n$ is even then the $\A_n(I)$ give a conformal net on $\Re$ that extends
to a M\"obius covariant positive energy pre-cosheaf on $S^1$.
 $\A_n$ satisfies duality on $S^1$. However $\A_n$ satisfies duality on $\Re$
{\it iff} $n=2$.
\par
 If $n$ is odd, $\A_n$ extends to a M\"obius covariant positive energy
pre-cosheaf on the disjoint union $S^1\sqcup S^1$. If $I$ is an interval in one
of the copies of $S^1$, we denote by $\tilde I$ the corresponding interval in
the other copy of $S^1$.
 Then duality holds in the form
 $$\A_n(I)'=\A_n(\tilde {I'})$$
 where $I'$ is the interior of the complement of $I$.
 \par
 Finally we notice that, with the previous notations, the tensor product
 $$\A_4\otimes\A_4$$
 is a conformally covariant  pre-cosheaf on $\Re^2$ such that Haag duality does
not hold.
 \section 3. Further results for Poincar\'e covariant
theories
 \par
 In this section we consider Poincar\'e covariant Quantum
Field Theories in  $d>2$ dimensions, namely a pre-cosheaf of
von~Neumann-algebras $\O\to\R(\O)$ where $\O$ is in the
family $\K$ of double cones in the Minkowski space $M$ with
the usual causality and additivity properties (see Section
1).
 \par
 The quasi-local $C^*$-algebra $\A_0$ generated by the local
algebras is supposed to act in the vacuum representation as usual
[\ref(Haag1)].
 \par
 The Poincar\'e group $\Poi$ acts by covariant automorphisms on $\A_0$ and
positivity of the energy-momentum is required.
 \par
 The pre-cosheaf is extended by additivity to general open regions in $M$. The
weak closure of the algebra associated with an unbounded region $\O$ will be
denoted by $\R(\O)$.
 \par
 We shall make two main assumptions in this Section:
 \item {$(a)$} Given any wedge region $W$, the modular
unitaries $\D_W^{it}$ of  $\R(W)$ act geometrically on the
pre-cosheaf $\A$,
$$\D^{it}_W\R(\O)\D^{-it}_W=\R(\L_t^W\O),
\qquad\O\in\K,\quad\forall W\in\W\eqno(3.1)$$
where $\L^W_t$ is the one-parameter group of diffeomorphisms defined in 1.6,
and $\W$ is the family of all wedge regions in $M$.
 \item{$(b)$} Distal split property holds, i.e there exist two regions
$\O_1\subset \O'_2$ in $M$ such that $\R(\O_1)$ and $\R(\O_2)$ generates a
$W^*$-tensor product.
 \par
 \claim{3.1 Theorem} If distal split property holds,
there is only one covariant unitary representation of the Poincar\'e group
on $\H$ leaving the vacuum vector invariant.
 \par
 \proof Distal split property implies that the group of internal symmetries
$\G$
is compact and commutes with any action of the Poincar\'e group by
automorphisms
[\ref(DoLo1)]. Then if two unitary covariant representations $U$, $V$
of $\Poi$ exist,  $\ad U(g)V(g^{-1})$ belong to $\G$, and therefore
gives rise to an action of $\Poi$ in $\G$. Since $\Poi$ has no non-trivial
finite-dimensional representations, and $\G$ is compact, such a
representation is trivial, i.e.  $\ad U(g)V(g^{-1})=id$. Then $U(g)V(g^{-1})$
is
a one-dimensional representation of $\Poi$, and, repeating the preceding
argument, $U(g)V(g^{-1})=1$.
 \par
 \claim{3.2 Lemma} Let $W$ be a wedge
in $\W$, $\D_W$ the modular operator for $(\R(W),\Q)$
and $U_W(t)$ the unitary associated with  $\L^W_t$ . Then, if assumptions $(a)$
and $(b)$ hold
$$z(t)=\D_W^{it}U_W(-t)$$
is a one-parameter group which commutes with the unitary
representation of $\Poi$, implements internal symmetries
and does not depend on the wedge $W$.
 \par
 \proof By assumption $(a)$, the unitaries
 $$z_t^W=\D^{-it}_WU_W(t)\quad,\qquad W\in\W,\quad t\in\Re$$
 implement internal symmetries, and therefore commute with the modular groups
of any local algebra.
 \np
 Moreover, if we fix $W_0\in\W$ and $W_1$ is any other wedge,
we can find an element $g\in\Poi$ which maps $W_1$ onto
$W_0$, and therefore
 $$z^{W_1}(t)=U(g^{-1})z^{W_0}(t)U(g)\qquad t\in\Re$$
 Since assumption $(b)$ implies that  the internal symmetries commute with any
action of $\Poi$ by automorphisms [\ref(DoLo1)], the unitary
$z^{W_1}(t)z^{W_0}(-t)\in\Co$. Since it preserves the vacuum vector, then
$z_t\equiv z_t^{W_0}$ is independent from $W_0$.
\np
Finally, the group property of $\{z(t)$, $t\in\Re\}$ is easily checked
using the mentioned commutation relations.
 \par
\claim{3.3 Theorem} Let $\A$ be a local pre-cosheaf on $M$
satisfying assumption $(a)$. Then the following holds:
\item{$(i)$} The pre-cosheaf $\A$ satisfies essential duality, i.e.
 $$\R(W)'=\R(W')$$
for each wedge region $W$.
\item{$(ii)$} If distal split property holds and $d>2$,
then
 $$\D_W^{it}=U_W(t)$$
where $\D_W$ is the modular operator for $(\R(W),\Q)$
and $U_W(t)$ is the  unitary representation of the boosts which
preserve $W$. The same result holds when $d=2$, provided that the parity
transformation is unitarily implemented.
 \par
 \proof The proof of part $(i)$ is identical to the proof of Theorem 2.3,
part~$(i)$.
 \np
 Since Lemma 3.4 holds, the proof of $(ii)$ is analogous to  the proof of
Theorem 2.3, part~$(ii)$, provided that  $R_i$ is replaced by a Poincar\'e
transformation mapping  a wedge $W$ onto $W'$. If $d>2$ we may
choose a rotation, if $d=2$ we use the parity transformation.
 \par
We conclude this section with an example of a non-split net where $(ii)$ of
Theorem 3.3 does not hold, namely $z(t)$ is non trivial.
\par
Such a construction was mentioned to us by D. Buchholz in a different context.
Let $\A$ be a Poincar\'e covariant net on the $d+1$-dimensional Minkowski space
$M_{d+1}$ which satisfies the Bisognano-Wichmann theorem. Define a net $\B$ on
the $d$-dimensional Minkowski space $M_d$ by
$$\B(\O)=\A(\p^{-1}\O)\qquad\O\subset M_d$$
where $\p:M_{d+1}\to M_d$ is the projection parallel to the $(d+1)$-th
coordinate.
\par
Let $W_1$ be a wedge in $M_d$, $z(t)$ the translation along the $(d+1)$ axis.
Then $z(t)$ implements a one-parameter group of internal symmetries of $\B$
which commutes with the action of the Poincar\'e group on $M_d$, hence
$$U'_W(t)\equiv z(t)U_W(t)$$
determines a new $d$-dimensional representation of $\Poi$ that violates the
Bisognano-Wichmann theorem.

\references{references4-dim}
\end